\documentclass{jnmp01b}

\usepackage{amsmath}

\numberwithin{equation}{section}

\newcommand {\Cee}    {{\mathbb  C}}

\newcommand {\fg}     {{\mathfrak{g}}}    %
\newcommand {\fa}     {{\mathfrak{a}}}
\newcommand {\fpo}    {{\mathfrak{po}}}
\newcommand {\fk}     {{\mathfrak{k}}}
\newcommand {\fsh}    {{\mathfrak{sh}}}
\newcommand {\fsl}    {{\mathfrak{sl}}}
\newcommand {\fgl}    {{\mathfrak{gl}}}  %
\newcommand {\fpsl}   {{\mathfrak{psl}}}
\newcommand {\fpsq}   {{\mathfrak{psq}}}
\newcommand {\fp}    {{\mathfrak{p}}}   %
\newcommand {\fpe}    {{\mathfrak{pe}}}   %
\newcommand {\fspe}   {{\mathfrak{spe}}}
\newcommand {\fas}    {{\mathfrak{as}}}
\newcommand {\fpsh}   {{\mathfrak{psh}}}
\newcommand {\fG}     {{\mathfrak{G}}}    %
\newcommand {\fA}     {{\mathfrak{A}}}

\newcommand {\Zee}    {{\mathbb  Z}}

\newcommand {\cC}     {{\cal C}}
\newcommand {\cM}     {{\cal M}}

\newcommand {\fq}     {{\mathfrak{q}}}     

\def \opname#1#2%
  {\expandafter\newcommand \csname #1\endcsname {{\mathop{#2}\nolimits}}}
\newcommand{\rmname}[1]
  {\expandafter\newcommand \csname #1\endcsname {{\operatorname{#1}}}}

\rmname{ad}
\rmname{tr}
\rmname{id}
\rmname{str}
\rmname{qtr}
\rmname{Res}

\opname{vvol}  {{v\hspace{-0.1ex}o\hspace{-0.02ex}l\/}}
\opname{Span} {{Span}}

\newcommand {\tto} {\longrightarrow}

\newtheorem{Statement}{Statement}[section]

\begin{document}

\renewcommand{\evenhead}{D.\ Leites and A.\ Shapovalov}
\renewcommand{\oddhead}{Manin-Olshansky Triples for Lie Superalgebras}


\thispagestyle{empty}

\begin{flushleft}
\footnotesize \sf
Journal of Nonlinear Mathematical Physics \qquad 2000, V.7, N~2,
\pageref{firstpage}--\pageref{lastpage}.
\hfill {\sc Letter}
\end{flushleft}

\vspace{-5mm}

\copyrightnote{2000}{D.\ Leites and A.\ Shapovalov}

\Name{Manin-Olshansky Triples for Lie Superalgebras}

\label{firstpage}

\Author{Dimitry LEITES and Alexander SHAPOVALOV}

\Adress{Deptartment of Mathathematics, University of Stockholm,
Roslagsv.\ 101,\\
Kr\"aftriket hus 6, 
SE-106 91, Stockholm, Sweden\\ 
E-mail: mleites@matematik.su.se}

\Date{Received September 18, 1999; Revised November 21, 1999; 
Accepted February 22, 2000}

\begin{abstract}
\noindent
Following V.~Drinfeld and G.~Olshansky, we construct
Manin triples $(\fg, \fa, \fa^*)$ such that $\fg$ is different from
Drinfeld's doubles of $\fa$ for several series of Lie superalgebras
$\fa$ which have no even invariant bilinear form (periplectic, Poisson
and contact) and for a remarkable exception.  Straightforward
superization of suitable Etingof--Kazhdan's results guarantee then the
uniqueness of $q$-quantization of our Lie bialgebras.
Our examples give solutions to the quantum Yang-Baxter equation in the
cases when the classical YB equation has no solutions.
To find explicit solutions is a separate (open) problem.  It is also
an open problem to list (\`a la Belavin-Drinfeld) all solutions of the
{\it classical} YB equation for the Poisson superalgebras $\fpo(0|2n)$
and the exceptional Lie superalgebra $\fk(1|6)$ which has a
Killing-like supersymmetric bilinear form but no Cartan matrix.
\end{abstract}


\section{Introduction}

Drinfeld proved \cite{D} that $\fa$ is a Lie bialgebra if and only if
it constitutes, together with its dual $\fa^*$, a {\it Manin triple}
$(\fg, \fa, \fa^*)$, where $\fg\cong\fa\oplus\fa^*$ as vector spaces,
$\fg$ is a Lie algebra possessing a nondegenerate invariant symmetric
bilinear form $B$ such that $\fa$ and $\fa^*$ are isotropic with
respect to $B$.  The known examples of forms $B$ arise from the Cartan
matrix of $\fg$ in case the Cartan matrix is symmetrizable.  Such a
form $B$ corresponds to the quadratic Casimir element $\Delta$ (which
for infinite dimensional $\fg$ belongs, strictly speaking, to a
completion of $U(\fg)$, rather than to $U(\fg)$ itself).  Etingof and
Kazhdan showed \cite{EK} how to $q$-quantize the bialgebra structure
of $U(\fa)$ in terms of $\fg$ and $\Delta$ and proved that such a
quantization is unique.

There are two ways to superize these results and constructions.  One
is absolutely straightforward.  Do not misread us!  We only refer to
the above-mentioned examples and result from deep and difficult papers
\cite{D} and \cite{EK}.  (Besides, performing the actual
implementation of these ``straightforward'' generalizations is quite a
job and, moreover, we encounter several unexpected phenomena, e.g.,
while presenting simple Lie superalgebras, cf.  \cite{GL2}.)

The only totally new feature of ``straightforward'' superization is
the fact that there is just one (in the class of $\Zee$-graded Lie
superalgebras of polynomial growth) series of simple Lie superalgebras
$\fsh(0|n)$ (see \cite{LS}) and one ``exception'', $\fk^L(1|6)$, (see
\cite{GL1}) which have no Cartan matrix but, nevertheless, possess an
invariant nondegenerate supersymmetric even bilinear form.

Another superization, performed by G.~Olshansky, is totally new. 
Drinfeld observed that the definition of Manin triples does not
require a nondegenerate symmetric invariant bilinear form on $\fa$. 
For a simple Lie algebra $\fa$, be it finite dimensional or of
polynomial growth, such a form always exists; hence, $\fa\cong\fa^*$
and $\fg\cong\fa\oplus\fa$, i.e., is the {\it double} of the bialgebra
$\fa$.  G.~Olshansky \cite{O} described two series of simple Lie
superalgebras $\fa$ without nondegenerate symmetric invariant bilinear
form for which there still exist Manin triples.  His construction
gives solutions to the {\it quantum} Yang-Baxter equation in the cases
the {\it classical} equation has {\it no} solutions, compare with the
classification \cite{LS}.

{\bf Our result}.  We list other examples of Manin triples of
G.~Olshansky type for the known Lie superalgebras of polynomial
growth.  Their $q$-quantization is routine thanks to
Etingof--Kazhdan's recipe.  (To get explicit formulas for the
$R$-matrix is, though routine, a quite tedious job and we leave it, so
far, as an open problem.)  Observe that in some of our examples, as
well as in \cite{O}, $\fa$ is not simple but \lq\lq close" to a simple
one (like affine Kac--Moody algebras).

\subsection{On background} 

For background on linear algebra in
superspaces and the list of simple stringy and vectorial Lie
superalgebras see \cite{Sh1}, \cite{Sh2}.  We recall here only the
additional data.  On a simple finite dimensional Lie algebra $\fg$
there exists only one up to proportionality invariant (with respect to
the adjoint action) nondegenerate symmetric bilinear form.  In this
one-dimensional space of invariant forms one usually chooses for the
point of reference the {\it Killing form} $ \langle x,
y\rangle_\ad=\tr(\ad x\cdot\ad y)$.  For any irreducible finite
dimensional representation $\rho$ of $\fg$ we have a proportional form
$\langle x, y\rangle_\rho=\tr(\rho(x)\cdot\rho(y))$.

Remark.  In reality, the Killing form is not the easiest to use, but
since all of the forms on the simple Lie superalgebra (of polynomial
growth or finite dimensional) are proportional to each other, one can
take most convenient for the task.  For $\fsl(n)$, for example, it is
more convenient to take the form associated with the standard
representation: $\langle x, y\rangle_\id=\tr(x\cdot y)$.

\begin{Statement} {\em (\cite{Sh2})} An invariant (with respect to the 
adjoint action) nondegenerate symmetric bilinear form on a simple Lie 
superalgebra $\fg$, if exists, is unique up to proportionality.
\end{Statement}

Observe that on Lie {\it super}algebras the form can be either 
even or odd, the Killing form can be identically zero even if there 
is another, nondegenerate, form.

Recall that for the matrices in the standard format 
$X=\begin{pmatrix}A&B\cr C&D\end{pmatrix}$ with elements in a 
supercommutative superalgebra $\cC$ the {\it supertrace} is defined as 
\[
\str: X=\tr\, A-(-1)^{p(X)}\tr\, D.
\]
Let $\fq(n; \cC)$ be the Lie superalgebra of 
supermatrices $Y=\begin{pmatrix}A&B\cr 
(-1)^{p(Y)+1}B&(-1)^{p(Y)}A\end{pmatrix}$, where $A, B\in \fgl(p|q; 
\cC)$ with $p+q=n$.  On it, $\str$ vanishes identically; instead, 
another invariant function, the {\it queertrace} is defined: 
\[
\qtr Y:=\str\, B.
\]
(On matrices in the standard format, i.e., if $pq=0$, 
the queertrace becomes $\qtr Y:=\tr\, B$.)

For the matrix Lie superalgebras the standard (identity) 
representation provides with nondegenerate forms given by the formula
\begin{equation}
\begin{split}
&\langle x, y\rangle_\id=\str(x\cdot y) \text{ is nondegenerate on
}\fgl(m|n);\\
&\langle x, y\rangle_\id=\qtr(x\cdot y) \text{ is nondegenerate on
}\fq(n)
\end{split}\label{(1.1.1)}
\end{equation}
and the forms induced by the above forms on the quotient algebras $\fpsl$ and
$\fpsq$ are also nondegenerate. 

Under the contraction $\fgl(2^n|2^n)\tto \fpo(0|2(n+1))$ and its 
restriction
$\fq(2^n)\tto \fpo(0|2n+1)$ (differential operator $\mapsto$ its 
symbol; we consider the Poisson bracket on the symbols) the traces 
(\ref{(1.1.1)}) turn into the integral and formulas (\ref{(1.1.1)}) turn into
\begin{equation}
\langle x, y\rangle_\id=\int(x\cdot y)~\vvol \text{ is nondegenerate on
}\fpo(0|N),\label{(1.1.2)}
\end{equation}
where $\vvol$ is the volume element (for its (nontrivial in
supersetting) definition see, e.g., \cite{L}). The form (\ref{(1.1.2)})  is even or odd
together with $N$; it induces a nondegenerate form on the simple Lie
superalgebra $\fsh(0|N)$.

\subsection{Bilinear forms on stringy superalgebras} 

Let $\fg$ be the 
Lie superalgebra $\fk(\cM^{2n+1|2n+6})$ of contact vector fields on 
$\cM$, where we consider either a compact supermanifold $\cM$ or 
contact vector fields generated by functions with compact support.  It 
so happens that for the dimension indicated $\fk(\cM^{2n+1|2n+6})$ 
preserves the volume element $\vvol$, see \cite{CLL}.  Define the form 
$B$ on $\fg$ setting $B(K_f, K_g)=\int fg~\vvol$.

In particular, consider the Fourier images of the functions on the 
supercircle $S^{1|6}$; let $\fk^L (1|6)$ be the Lie superalgebra of 
the corresponding contact fields.  It differs from $\fk(1|6)$ by the 
possibility to consider negative powers of $t=\exp i\varphi$, where 
$\varphi$ is the angle parameter on $S$; the superscript $L$ indicates 
that we consider Laurent coefficients, rather than only polynomial ones.  
We set
\[
B(K_f, K_g)=\Res fg,\;  \text{ where }\Res ~(f) = \text{the coefficient of
}\; \frac{\xi_1\xi_2\xi_3\eta_1\eta_2\eta_3}{t}. 
\]
In \cite{GL1} the even quadratic Casimir element $\Delta$ for $\fk^L 
(1|6)$ corresponding to $B$ is explicitly computed.  In terms of this 
element a solution of the classical Yang-Baxter equation can be 
expressed in the same way as in \cite{LS}, namely, 
$\frac{\Delta}{u-v}$, and from this solution a quantum solution can be 
uniquely recovered thanks to the general uniqueness theorems of 
\cite{CP}, \cite{EK}.

In \S 2 we offer other solutions, without classical counterpart.

\subsection{Open problems} 

It remains to explicitly produce the 
universal $R$-matrix for the triples of \S 2 (sec. 2.2 --2.6)
and classify the 
trigonometric solutions of the classical YB equation with values in  
$\fk^L (1|6)$ and in the Poisson superalgebras $\fpo(0|2n)$, a problem 
left open in \cite{LS}.

\section{Main result: examples of Manin triples for 
Lie superalgebras}

{\bf 2.0}.  It is known that the extension of Drinfeld's Example 3.2 
in \cite{D} to a simple Lie {\it super}algebra $\fa$ with a 
symmetrizable Cartan matrix (or, equivalently, an even supersymmetric 
nondegenerate form) has a new feature: there are {\it several} 
nonisomorphic Borel subsuperalgebras $\fa^*$; but, nevertheless, 
$\fg\cong\fa\oplus \fa^*$ for all these Borel subalgebras.

(Observe that this implies that we can take for $\fa$ nonisomorphic 
Borel subalgebras and get various types of defining relations for 
$\fg$.  This observation applies not only to the superization of 
Example 3.2 in \cite{D} --- finite dimensional Lie superalgebras with 
a symmetrizable Cartan matrix --- but to infinite dimensional 
superizations and the (analogues of) Borel subsuperalgebras of 
$\fpo(0|2n)$ and $\fk^L(1|6)$ considered below.  For a classification 
of systems of simple roots and for a description of Borel 
superalgebras see \cite{PS}, for a description of defining relations 
see \cite{GL2} and refs.  therein.)

In their famous papers \cite{KT} Khoroshkin and 
Tolstoy {\it explicitly} wrote the universal $R$-matrices.

{\bf 2.1}.  (\cite{O}) $\fa=\fq(n)$; let $\fa^*=\left \{\pmatrix 
X&Y\\0&Z\endpmatrix\right\}$, where $X$ is upper triangular, $Z$ is 
lower triangular and $X_{ii}=-Z_{ii}$; then $\fg=\fgl(n|n)$ with the 
form $B(x, y)=\str (xy)$.  In this example the Lie superalgebra $\fa$ 
is not simple.  The construction can be transported to the simple 
algebra $\fa=\fpsq (n)$ (cf.  Remark 5 of \cite{O}):

{\bf 2.1}$'$.  $\fa=\fpsq(n)$; let $\fa^*$ be as in Example 2.1 but $\tr 
X=\tr Z=0$; then $\fg=\fp\fsl(n|n)$ with the form induced by $B(x, y)=\str 
xy$. (In the last formula we consider the conditional presentation of the 
elements from $\fp\fsl(n|n)$ by supermatrices from $\fsl(n|n)$ and 
the modified bracket: composition of the bracket with the subsequent 
subtraction of the supertrace.)

{\bf 2.2}.  $\fa=\fpe(n)$; let $\fa^*=\left\{\pmatrix 
X&Y\\Z&T\endpmatrix\right\}$, where $X$, $Z$ and $T$ are upper 
triangular matrices and $X_{ii}=T_{ii}$, $Y$ is a strictly upper 
triangular matrix; $\fg=\fgl(n|n)$ with the form $B(x, y)=\str xy$.

{\bf 2.2}$'$.  $\fa=\fspe(n)=\fpe(n)\cap\fsl(n|n)$; let $\fa^*$ be as 
in our Example 2.2 with $\tr X=\tr T=0$; then $\fg=\fp\fsl(n|n)$ with the 
form as in Example 2.$1'$.

{\bf 2.3}.  $\fa=\fpo(0|2n-1)$.  Let $\fg=\fpo(0|2n)$ be generated 
by $\Cee[\theta_1, \dots \theta_{2n}]$ and $\fa$ by $\Cee[\theta_1, 
\dots \theta_{2n-1}]$.  The form $B$ on $\fg$ is $B(x, y)=\int 
xy~\vvol(\theta)$.  Let $\fa^*=\fa\cdot \zeta$, where 
$\zeta=\theta_{2n}+\sum\limits_{1\leq i\leq 2n-1}k _i\theta_i$ and 
$\sum\limits_{1\leq i\leq 2n-1}k _i^2=-1$.

{\bf 2.3}$'$.  $\fa=\fpsh(2n-1)$.  Let $\fg=\fpsh(0|2n)$ be 
conditionally (as in 2.1$'$) realized by generating functions 
$\Cee[\theta_1, \dots \theta_{2n}]$ and $\fa$ by $\Cee[\theta_1, \dots 
\theta_{2n-1}]$.  The form $B$ on $\fg$ is the one induced by the 
Berezin integral.  Then $\fa^*=\fa \cdot \zeta$, with the 
$\zeta$ from 2.3.

{\bf 2.4}.  $\fg=\fpo(0|6)$.  Let $\fa=\fas=\Span(1, \Lambda(\xi, 
\eta), \Lambda^2(\xi, \eta), g_{11}^\xi)$, see \cite{Sh1}; 
then $\fa^*$ is generated (as Lie superalgebra) by $g_{11}^\eta$
(see [12]) and 
$\xi_1\xi_2\xi_3\eta_1\eta_2\eta_3$.

{\bf 2.5}. Manin triples for $\fg=\fk^L(1|6)$.
(Recall that superscript $L$ indicates that we consider Laurent
polynomials as generating functions, not just polynomials as in
2.5.1.)

{\bf 2.5.1}.  $\fa$ and $\fa^*$ as in Drinfeld's Example 3.3 
\cite{D}, i.e., $\fa=\fk(1|6)\cong\Span(K_f: f\in \Cee[t, \xi, \eta])$ 
and $\fa^*=\Span(K_f: f\in t^{-1}\cdot\Cee[t^{-1}, \xi, \eta])$.

{\bf 2.5.2}. $\fa$ and $\fa^*$ as in Example 3.2 \cite{D}, for various
Borel subalgebras. 

{\bf 2.5.3}.  A variation of our Example 2.3: $\fa=\fk^L(1|5)$ and 
$\fa^*=\fa\cdot \zeta$, where $\zeta=\theta_{6}+\sum\limits_{1\leq 
i\leq 5}k_i\theta_i$ and $\sum\limits_{1\leq i\leq 5}k _i^2=-1$.

{\bf 2.6}.  In what follows $V^{1}=V\otimes \Cee[x^{-1}, x]$ for any 
vector space $V$.  The form $B^{(1)}$ on $\fg^{(1)}$ is $\Res_x 
B(f(x), g(x))$, where $x=\exp i\psi$ for the angle parameter $\psi$ on 
$T^1$ and $B$ is the nondegenerate invariant supersymmetric form on 
$\fg$.  (Clearly, this $T^1$ and $S$ from sec.  1.2 are diffeomorphic 
circles, as any two circles are.  We denote them differently to 
underline that these are {\it different} circles.)  Then new Manin 
triples $(\fG, \fA, \fA^{*})$ are obtained from the triples $(\fg, 
\fa, \fa^{*})$ corresponding to the examples 2.1--2.4 above for 
$\fG=\fg^{(1)}$, $\fA=\fa^{(1)}$, $\fA^{*}=(\fa^{*})^{(1)}$.

\subsection*{Acknowledgements}

We most gratefully acknowledge financial support of NFR and
thank G.~Olshansky for help.

\label{lastpage}

\end{document}